\documentclass{article}

\usepackage{arxiv}
\usepackage[utf8]{inputenc} 
\usepackage[T1]{fontenc}    
\usepackage{hyperref}       
\usepackage{url}            
\usepackage{booktabs}       
\usepackage{amsfonts}       
\usepackage{nicefrac}       
\usepackage{microtype}      
\usepackage{lipsum}
\usepackage{amsmath,amssymb}
\usepackage{amsthm}
\usepackage{graphicx}
\usepackage{listings}
\usepackage{color}
\usepackage{fancyhdr}
\usepackage[framemethod=tikz]{mdframed}
\usepackage{subfig}
\usepackage[nottoc,notlot,notlof]{tocbibind}
\usepackage{csvsimple}
\usepackage{longtable}
\usepackage{colortbl}
\usepackage{float}
\usepackage[linesnumbered,lined,boxed,commentsnumbered]{algorithm2e}
\usepackage{verbatim}
\usepackage{lscape}
\usepackage{multirow}
\usepackage{hhline}

\newcommand{\R}{\mathbb{R}}
\newcommand{\N}{\mathbb{N}}

\DeclareMathOperator{\capop}{cap}

\def\restriction#1#2{\mathchoice
              {\setbox1\hbox{${\displaystyle #1}_{\scriptstyle #2}$}
              \restrictionaux{#1}{#2}}
              {\setbox1\hbox{${\textstyle #1}_{\scriptstyle #2}$}
              \restrictionaux{#1}{#2}}
              {\setbox1\hbox{${\scriptstyle #1}_{\scriptscriptstyle #2}$}
              \restrictionaux{#1}{#2}}
              {\setbox1\hbox{${\scriptscriptstyle #1}_{\scriptscriptstyle #2}$}
              \restrictionaux{#1}{#2}}}
\def\restrictionaux#1#2{{#1\,\smash{\vrule height .8\ht1 depth .85\dp1}}_{\,#2}}

\mdfdefinestyle{thm}{
	nobreak=true,
	hidealllines=true,
	leftline=true,
	innerleftmargin=10pt,
	innerrightmargin=10pt,
	innertopmargin=0pt,
}

\surroundwithmdframed[style=thm]{definition}
\newtheorem{proposition}{Proposition}[section]\surroundwithmdframed[style=thm]{proposition}
\newtheorem{theorem}{Theorem}[section]\surroundwithmdframed[style=thm]{theorem}
\surroundwithmdframed[style=thm]{lemma}
\surroundwithmdframed[style=thm]{property}
\newtheorem{problem}{Problem}[section]\surroundwithmdframed[style=thm]{problem}
\newtheorem*{problem*}{Problem}\surroundwithmdframed[style=thm]{problem*}

\theoremstyle{remark}

\title{Iterative Approach to Image Compression with Noise : Optimizing Spatial and Tonal Data}

\author{
  Zakaria BELHACHMI \\ 
  IRIMAS\\
  Université de Haute-Alsace\\
  Mulhouse, France \\
  \texttt{zakaria.belhachmi@uha.fr} \\
   \And
  Thomas JACUMIN \\
  IRIMAS \\
  Université de Haute-Alsace\\
  Mulhouse, France \\
  \texttt{thomas.jacumin@uha.fr} \\
}

\begin{document}
\maketitle

\begin{abstract}
We consider some iterative methods for finding the best interpolation data in the images compression with noise. The interpolation data consists of the set of pixels and their grey/color values. The aim in the iterative approach is to allow the change of the data dynamically during the inpainting process for a reconstruction of the image that includes the enhancement and denoising effects. The governing PDE model of this approach is a fully parabolic problem where the set of stored pixels is time dependent. We consider the semi-discrete dynamical system associated to the model which gives rise to an iterative method where the stored data are modified during the iterations for best outcomes. Finding the compression sets follows from a shape-based analysis within the $\Gamma$-convergence tools developed in \cite{Belhachmi2009, Belhachmi2022}, in particular well suited topological asymptotic and a ``fat pixels'' approach are considered to obtain an analytic characterization of the optimal sets in the sense of shape optimization theory. We perform the analysis and derive several iterative algorithms that we implement and compare to obtain the most efficient strategies of compression and inpainting for noisy images. Some numerical computations are presented to confirm the theoretical findings. 
Finally, we propose a modified model that allows the inpainting data to change with the iteration and compare the resulting new method to the ``probabilistic'' ones from the state-of-the-art.
\end{abstract}


\section*{Introduction}

Image compression and \textit{inpainting} aim to reconstruct missing parts of an image from a chosen set of few known pixels, called \textit{inpainting mask} \cite{Masnou1998, Bertalmio2000}.
PDE-based models have emerged among the most efficient methods in this field as they allow high quality reconstruction even from small set of pixels (masks) \cite{Schmaltz2009,Galic2005} and the references therein.  
Actually, diffusion-based models \cite{Galic2008,Schmaltz2009,Belhachmi2009} are known to be efficient in image processing, particularly in image and video compression \cite{Andris2016, Andris2021}, and similar problems such as zooming or denoising \cite{Adam2017, Lenzen2011}.

It is commonly admitted in image compression problems, that both a well chosen inpainting operator and mask  play a crucial role on the quality of the reconstruction and the efficiency of the compression strategy \cite{Mainberger2012,Bae2010, Galic2008}. In fact, one aims to use a simple differential operator (e.g. laplacian) for easy and cheap reconstruction, but without sacrificing important features of the image at hand (e.g. singularities -edges-). This may appears somehow contradictory (because of the regularizing effects of the operator) or at least calls for finding a good balance between the mentioned simplicity of the differential operator and the necessary care about too strong regularization. Now, fixing the operator, the main question is the following : is a good ``choice'' of the mask possible? and what is a good choice? Many studies show that apart from noise or textures, good mask candidates should include pixels from contours of objects in the image. A mathematical theory addressing such questions from the shape optimization point of view may be found in \cite{Belhachmi2009}, adapted to noisy images in \cite{Belhachmi2022}. In particular, it is proven there that an optimal mask under  ``size constraints'' exists and is easily characterized by an analytical criterium. The most significant advantage then being the selection of the mask for general meshes and without any cost (other than the storage of the pixels).
%
%

Recently, many research works (\cite{Schmaltz2014,Chizhov2021,Belhachmi2022}) noticed that choosing the mask at once may be less efficient than improving it somehow iteratively, however without giving a clear and sound justified strategy for the iterative process.

In this article, we consider the linear heat equation as the inpainting operator and study the selection of a ``time-dependent'' set of optimal pixels using shape optimization tools, extending this way the results in \cite{Belhachmi2022}. We present several strategies combining the shape optimization approach giving analytic criteria for masks selection and an enhancement process allowing to enrich/adapt such choice of pixels as well as the modification of the data on this points (\textit{tonal optimization}). This turns to be a nice tool for image compression, particularly those with noise.

\subsection*{Related works}

There exists important literature and several works addressing the problem of images compression (see \cite{Distasi1997, Dell2006, Mainberger2012, Hoeltgen2013, Schmaltz2014} and the references therein). Most of these methods require a significant amount of computation time for the masks selection as the set of pixels selected are built upon a process based on local choices at the scale of each pixel. Moreover, mostly, the notion of optimality of the selection criteria, in the mathematical sense, is barely considered, though some of these methods give very good results and may include powerful features such as a low storage cost, handling complex images (e.g. with noise or textures). 
%
Beside finding an optimal mask, \cite{Mainberger2012} introduced a data optimization strategy allowing to modify the values stored on the inpainting mask to improve the reconstructed images. Loosely speaking, they optimize the spatial distribution of the inpainting data, with a probabilistic data sparsification followed by a nonlocal pixel exchange (randomly remove and add pixels in the mask to avoid being trapped in a local minimizer) and then they optimize the grey values in these inpainting points using a least squares approach. Other data optimization approaches have been proposed in \cite{Chizhov2021} by using finite element method instead of finite differences method and deep learning techniques for both choosing the inpainting mask and values \cite{Alt2022}.

In this article, we aim to propose some iterative methods using both the mathematical results in the spirit of \cite{Belhachmi2009} and an iterative process to enhance the quality of the mask within a data optimization approach for images with noise. We study and analyze the problem of finding a set $K$ for the time harmonic linear diffusion that is ``time-dependent'', in the sense of allowing mask's changes (by exchanging, removing, or adding pixels) and changes of the stored values during the process. We propose a deterministic  all-in-one strategy for the mask choice and the tonal optimization. The enhancement of the mask selection is encoded in the criterium which is dependent on the previous values of the reconstructed image and the tonal optimization is performed by replacing the original values of the noisy image with those given by the iterative reconstruction. 
%
We compare different methods proposed here and those existing in the related literature, particularly in the presence of noise. 

The standard PDE-based compression problem reads : let $D\subset\R^2$ the support of an image (say a rectangle) and $f : D\longrightarrow \R$, be an image which is assumed to be known only on some region $K\subset D$. There are several PDE models to interpolate $f$ in $D\setminus K$, the simplest being the use of the linear elliptic equation, with Dirichlet boundary data $f\vert_{K}$ on $K$ and homogeneous Neumann boundary conditions on $\partial D$. In our approach, we consider the iterative system of equations 

\begin{problem} For $n\in\N$, given $u^n$, find $u^{n+1}$ in $H^1(D)$ such that
	\begin{equation}
		\left\{\begin{array}{rl}
			u^{n+1} - \alpha \Delta u^{n+1} = u^n, & \text{in}\ D\setminus K_n, \\
			u^{n+1} = f, & \text{in}\ K_n, \\
			\frac{\partial u^{n+1}}{\partial \mathbf{n}} = 0, & \text{on}\ \partial D, \\
		\end{array}\right .\label{eq:l2insta}
	\end{equation}
	\label{pb:l2insta}
\end{problem}
with the initial conditions $u_0$ and $K_0$. The data $f$ in $K_n$ may be replaced during the iteration by a function $h_n$, possibly depending on $n$, for tonal optimization and which in our case might be a ``smoothed'' version of $f$ upon the iterations.
Formally, the family of Problems \ref{pb:l2insta}, with $\alpha=\delta t$, the time step, is a semi-implicit discrete system to solve \\

\begin{problem} For $t>0$, find $u(t,\cdot)$ in $H^1(D)$ such that
	\begin{equation}
		\left\{\begin{array}{rl}
			\partial_t u(t,\cdot) - \Delta u(t,\cdot) = 0, & \text{in}\ D\setminus K_t, \\
			u(t,\cdot) = f, & \text{in}\ K_t, \\
			\frac{\partial u(t,\cdot)}{\partial \mathbf{n}} = 0, & \text{on}\ \partial D, \\
		\end{array}\right .\label{eq:linear_diffusion_filter}
	\end{equation}
	\[ u(0,\cdot) = u_0,\ \text{in}\ D, \]
	\label{pb:linear_diffusion_filter}
\end{problem}

It might be possible and interesting to make precise the word ``formally'' by providing a well suited topology of sets and addressing the question of convergence of the semi-discrete system to Problem \ref{pb:linear_diffusion_filter} but this is a quite deep question beyond the scope of this article.

Problem \ref{pb:l2insta} with a fixed $K$ has been considered in \cite{Belhachmi2022} and may appears as  a particular case of the iterative process where one look for an optimal inpainting mask once and for all. The goal in solving Problem \ref{pb:l2insta} is to modify $K_n$ according to the change of the data (we notice that $u^n$ is less noisy than $f$) while preserving the optimality at each step $n$. Thus, the iterative approach leads to adaptive improvement of the whole data with : topological changes in $K_n$ and stored values changes that is $h_n\mid_{K_n}$ which may differ from $f$ (\textit{tonal optimization}). 

We assume given $f\in H^1(D)$ and $\Delta f\in L^2(D)$ with $\frac{\partial f}{ \partial\mathbf{n}}=0$,  for simplicity though in practice $f$ is a function of bounded variations with a non trivial jump set. In fact, the whole analysis in the paper extends to the case of $f\in L^2(D)$. 

To ensure the first compatibility condition for the full parabolic system, with the non-homogeneous ``boundary'' condition, we take $u(0,.)=f$ in $D$. Note that this choice is not in contradiction with image compression since the entire image is available during the encoding step. Denoting by $u_{K_n}$ the solution of Problem \ref{pb:l2insta}, the question is to identify the region $K_n$ which gives the ``best'' approximation $u_{K_n}$, in a suitable sense, for example which minimizes some $L^p$ or Sobolev norms, e.g. in \cite{Belhachmi2009}
\[ \int_D \vert\nabla u_{K_n}-\nabla f\vert^2\ dx, \]
For noisy images, a well suited  choice for the reconstruction is to minimize the $L^p$-norms of $u_{K_n}-f$, particularly for $p=1$ and $p=2$, which are known to be good filters for a large class of additive noises.
In this article, we will take $p=2$ for easy computations of the topological gradient but the main analysis holds for $p=1$. 

Following \cite{Belhachmi2009}, we develop at each step $n>0$, two methods of finding an optimal shape. The first is based on a topological asymptotic (a gradient method) which gives pointwise information on the set of pixels to select and yields a hard thresholding approach. The second method consists in a finite dimensional shape optimization problem, where $K_n$ is taken as a union of small balls, i.e. a finite number of ``fat pixels''. Then, performing the asymptotic analysis by $\Gamma$-convergence when the number of pixels is increasing (in the same time that the fatness vanishes), we obtain useful information about the optimal distribution of the best interpolation pixels as a density function leading to soft thresholding approach.

In all cases, we obtain an optimal mask, in the sense  mentioned above, and the values to be stored within the iterative method.

\subsection*{Organisation of the article}

In Section \ref{sec:l2insta:continuous_model}, we recall the mathematical model of the compression problem and we summarize the analysis steps and results referring interested readers to \cite{Belhachmi2009,Belhachmi2022} for details and recall the two methods proposed to construct our set of interpolation points at each step $n$. In Section \ref{sec:algorithms}, we propose three different algorithms based on the hard/soft-thresholding criterium stated in the previous section. In Section \ref{sec:comparison} we numerically compare the proposed algorithms for both image compression and image denoising. Finally, in Section \ref{sec:numerical_results}, we propose a modified model that allows the inpainting data to change with the iteration and compare the resulting new methods with the so-called \textit{sparsification} and \textit{densification} algorithm \cite{Adam2017,Mainberger2012}.

\section{Review of the Continuous Model}
\label{sec:l2insta:continuous_model}

In this section, we mainly recall and adapt results found in \cite{Belhachmi2022} in the case of the non-homogeneous linear diffusion inpainting.

\subsection{Min-max Formulation}

Let $D$ be a smooth bounded open subset of $\R^2$. The shape optimization problem we study is, for a given $n\in\N$,
\begin{equation}
	\min_{K_n\subseteq D,\ \capop(K_n)\leq c}\Big\{  \frac{1}{2}\int_D |u_{K_n}-f|^2\ dx + \frac{\alpha}{2}\int_D |\nabla (u_{K_n}-f)|^2\ dx\ \Big|\ u_{K_n}\ \text{solution of Problem}\ \ref{pb:l2insta} \Big\},
	\label{pb:l2insta_opt_no_constraint}
\end{equation}
where $\capop(E)$ is the capacity of a subset $E$ in $D$ i.e.
\[ \capop(E) = \inf\Big\{ \int_D |\nabla u|^2\ dx + \int_D u^2\ dx\ \Big|\ u\in H^1_0(D),\ u\geq 1 \text{ a.e. in } E \Big\}, \] and $c>0$. The image compression problem aims to find an optimal set of pixels from which an accurate reconstruction of the (noisy) image will be performed. If we denote by $u_{K_n}$ the solution of Problem \ref{pb:l2insta}, it is straightforward to obtain \\
\begin{proposition}
	\label{prop:l2insta_opt_no_constraint:reformulation}
	The optimization problem \eqref{pb:l2insta_opt_no_constraint} is equivalent to 
    \begin{equation}
    	\label{pb:l2insta_opt_penalized}
    	\max_{K_n\subseteq D,\ \capop(K_n)\leq c} \min_{u\in H^1(D),\ u=f\ \text{in}\ K_n} \frac{\alpha}{2}\int_D |\nabla u|^2\ dx + \frac{1}{2} \int_D (u-u^n)^2\ dx - \beta m(K_n),
    \end{equation}
    for $\beta>0$.
\end{proposition}
The next section is devoted to the analysis within the $\gamma$-convergence approach follows the same lines as in \cite{Belhachmi2009, Belhachmi2022} with slight changes.

\subsection{Analysis of the Model}

It is well known that such shape optimization problems do not always have a solution (e.g. \cite{Belhachmi2009}), we seek a relaxed formulation, which yields a relaxed solution, that is to say, a capacity measure. Thus, we consider the relaxed problem
\begin{equation}
    \max_{\mu_n\in\mathcal{M}_0(D)} \min_{u\in H^1(D)} \frac{\alpha}{2}\int_D |\nabla u|^2\ dx + \frac{1}{2} \int_D (u-u^n)^2\ dx + \frac{1}{2}\int_D (u-f)^2\ d\mu_n - \beta\capop(\mu_n),
    \label{eq:l2insta_opt_penalized_relaxed}
\end{equation}
where $\mathcal{M}_0(D)$ is the set of all non negative Borel measures $\mu$ on $D$, such that
\begin{itemize}
	\item $\mu(B)=0$, for every Borel set $B$ subset of $D$ with $\capop(B) = 0$,
	\item $\mu(B)=\inf\big\{\mu(U)\ \big|\ U\ \text{quasi-open},\ B\subseteq U\big\}$, for every Borel subset $B$ of $D$,
\end{itemize}
and $\capop(\mu)$ is the measure capacity i.e. for $\mu$ in $\mathcal{M}_0(D)$,
\[ \capop(\mu) := \inf_{u\in H^1_0(D)} \int_D |\nabla u|^2\ dx + \nu\int_D u^2\ dx + \int_D (u-1)^2\ d\mu. \]
Next we give a natural way to identify a set to a measure of $\mathcal{M}_0(D)$. Let $E$ be a Borel subset of $D$. We denote by $\infty_E$ the measure of $\mathcal{M}_0(D)$ defined by
\[ \infty_E(B) := \begin{cases}
	+\infty &,\ \text{if}\ \capop(B\cap E)>0,\\
	0 &, \ \text{otherwise}.
\end{cases},\ \text{for all}\ B\ \text{Borel subset of}\ D. \]
Formally, when $\mu_n = \infty_A$ in \eqref{eq:l2insta_opt_penalized_relaxed}, we penalize the functional which is equal to $+\infty$ until $u=f$ in $A$. Somehow, we embed the Dirichlet condition into the min-max problem thanks to this penalization. For technical reasons, we want to include balls centered at points $x_0$ in $D$, that we do not want to be too close to the boundary of $D$, we introduce the following notations for $\delta>0$,
\[ D^{-\delta} := \{ x\in D\ |\ d(x,\partial D) \geq \delta \}\subseteq D, \]
\[ \mathcal{K}^\delta(D) := \{ K\subseteq D\ |\ K\ \text{closed},\ K\subseteq D^{-\delta} \}, \]
and
\[ \mathcal{M}_0^\delta(D) := \{ \mu\in\mathcal{M}_0(D)\ |\ \restriction{\mu}{D\setminus D^{-\delta}}=0 \} \subseteq\mathcal{M}_0(D). \]
We have that \eqref{eq:l2insta_opt_penalized_relaxed} is the relaxed formulation of the optimization problem \eqref{pb:l2insta_opt_penalized} in the sense of the $\gamma$-convergence (Theorem 1.4 in \cite{Belhachmi2022}) : \\
\begin{theorem}
	We have,
	\[ \sup_{K\in\mathcal{K}_\delta(D)} \big( E(\infty_K) - \beta\capop(\infty_K) \big) = \max_{\mu\in\mathcal{M}_0^\delta(D)}\big( E(\mu) - \beta\capop(\mu) \big), \]
	where
	\[ E(\mu) := \min_{u\in H^1(D)} F_\mu(u), \]
	and
	\[ F_\mu(u) := \begin{cases}
    	\alpha\int_D |\nabla u|^2\ dx + \int_D (u-f)^2\ d\mu, &\text{if}\ |u| \leq |f|_\infty, \\
    	+\infty, &\text{otherwise.}
    \end{cases} \]
\end{theorem}
This result gives us the existence of a relaxed solution which $\gamma$-converges to the ``solution'' of our shape optimization problem. In order to solve the relaxed problem, we may use a shape derivative with respect to the measures $\mu_n$. However, such a method yields diffuse measures, thus too thick sets whereas we seek discrete sets of pixel. 

\subsection{Topological Gradient}

Here, we aim to compute the solution of our optimization problem \eqref{pb:l2insta_opt_no_constraint} by using a topological gradient-based algorithm as in \cite{Larnier2012, Garreau2001}. This kind of algorithm consists in starting with $K_n = \bar{D}$ and determining how making small holes in $K_n$ affect the cost functional to find the balls which have the most decreasing effect. To this end, let us define $K_\varepsilon^n$ the compact set $K_n\setminus B(x_0,\varepsilon)$ where $B(x_0,\varepsilon)$ is the ball centered in $x_0\in D$ with radius $\varepsilon>0$ such that $B(x_0,\varepsilon)\subset K_n$. From now, we consider the variable $v_{K_n} := u_{K_n} - f$, with $u_{K_n}$ solution of Problem \ref{pb:l2insta}. Let us denote by $j$ the functional \[ j : A\subset D \mapsto \min_{v\in H^1(D),\ v=0\ \text{in}\ A} \frac{\alpha}{2}\int_D |\nabla v|^2\ dx + \frac{1}{2} \int_D v^2\ dx\ - \int_D g^n v\ dx, \]
where $g^n:= u^n-f+\alpha\Delta f$. We denote by $v_\varepsilon^n$ the minimizer of $j(K_\varepsilon^n)$ and we have \\
\begin{proposition} With notations from above, we have when $\varepsilon$ tends to $0$,
	\[ j(K_\varepsilon^n) - j(K_n) = \frac{\pi}{2}\big(g^n(x_0)\big)^2\varepsilon^2\ln(\varepsilon) + O(\varepsilon^2). \]
	\label{prop:topologicalGradient}
\end{proposition}
Since for $\varepsilon < 1$, $\ln\varepsilon < 0$, the result above suggests to keep the points $x_0$ where $ |u^n(x_0)-f(x_0)+\alpha\Delta f(x_0)|^2$ is maximal, when $\varepsilon$ small enough. In the next section, we will see that such a strict threshold rule might be relaxed. 

\subsection{Optimal Distribution of Pixels : The ``Fat Pixels'' Approach}

In this section, we change our point of view by considering ``fat pixels'' instead of a general set of interpolation points. In the sequel, we will follow \cite{Buttazzo2005, Belhachmi2009, Belhachmi2022}. We restrict our class of admissible sets as an union of balls which represent pixels. For $m>0$ and $k\in\N$, we define
\[ \mathcal{A}_{m,k} := \Big\{ \overline{D}\cap\bigcup_{i=1}^k\overline{B(x_i,r)}\ \Big|\ x_i\in D_r,\ r=mk^{-1/2} \Big\}, \]
where $D_r$ is the $r$-neighborhood of $D$. We consider problem \eqref{pb:l2insta_opt_no_constraint} for every $K_n\in\mathcal{A}_{m,k}$ i.e. 
\begin{equation}
	\min_{K_n\in\mathcal{A}_{m,k}}\Big\{ \frac{1}{2} \int_D (u_{K_n}-f)^2\ dx + \frac{\alpha}{2}\int_D |\nabla u_{K_n}-\nabla f|^2\ dx\ \Big|\ u_{K_n}\ \text{solution of Problem}\ \ref{pb:l2insta} \Big\}.
\end{equation}
By setting $v_{K_n} := u_{K_n} - f$, this last optimization problem can be reformulated as a compliance optimization problem. We set $g^n := u^n - f + \alpha\Delta f$ like in the previous section. Here, we do not need to specify a size constraint on our admissible domains. Indeed, imposing $K_n\in\mathcal{A}_{m,k}$ implies a volume constraint and a geometrical constraint on $K_n$ since $K_n$ is formed by a finite number of balls with radius $m\,k^{-1/2}$. The well-posedness of such a problem has been studied in the laplacian case in \cite{Buttazzo2005}. As pointed out in \cite{Bucur2005}, the local density of $K_k^{n,\text{opt}}$, the optimal solution, can be obtained by using a different topology for the $\Gamma$-convergence of the rescaled energies. In this new frame, the minimizers are unchanged but their behavior is seen from a different point of view. We define the probability measure $\mu_K$ for a given set $K$ in $\mathcal{A}_{m,k}$ by
\[ \mu_K := \frac{1}{k}\sum_{i=1}^k \delta_{x_i}. \]
We define the functional $F_k^n$ from $\mathcal{P}(\bar D)$ into $[0,+\infty]$ by
\[ F_k^n(\mu) := \begin{cases} k\int_D g^n v_{K_n}\ dx &,\ \text{if}\ \exists K_n\in\mathcal{A}_{m,k},\ \text{s.t.}\ \mu=\mu_{K_n}, \\
+\infty&,\ \text{otherwise}. \end{cases} \]
The following $\Gamma$-convergence of $F_k^n$ theorem is similar to the one given in Theorem 2.2. in \cite{Buttazzo2005}. \\
\begin{theorem}
    \label{thm:g-convergence}
	Then the sequence of functionals $F_k^n$, defined above, $\Gamma$-converge when $k$ tends to $+\infty$ with respect to the weak $\star$ topology in $\mathcal{P}(\bar{D})$ to 
	\begin{equation}
	    F^n(\mu^n) := \int_D \frac{(g^n)^2}{\mu_a^n}\theta\big(m(\mu_a^n)^{1/2}\big)\ dx,
	    \label{eq:g-convergence}
	\end{equation}
	where $\mu^n = \mu_a^n dx + \nu^n$ is the Radon-Nikodym-Lebesgue decomposition of $\mu^n$ (\cite{Folland2013}, Theorem 3.8) with respect to the Lebesgue measure and
	\[ \theta(m) := \inf_{K_k\in\mathcal{A}_{m,k}} \liminf_{k\to +\infty} k \int_D g^n v_{K_k}\ dx, \]
	$v_{K_n} := u_{K_n} - f$, $u_{K_n}$ solution of Problem \ref{pb:l2insta}.
\end{theorem}
As a consequence of the $\Gamma$-convergence stated in the theorem above, the empirical measure $\mu_{K_k^{n,\text{opt}}} \to \mu^{n,\text{opt}}$ weak $\star$ in $\mathcal{P}(\R^d)$ where $\mu^{n,\text{opt}}$ is a minimizer of $F^n$. Unfortunately, the function $\theta$ is not known explicitly. Formal Euler-Lagrange equation and the estimates on $\theta$ \cite{Belhachmi2022} give the following information : to minimize \eqref{eq:g-convergence} one have to take
\[ \frac{(\mu_a^n)^2}{|1-\log\mu_a^n|} \approx c_{m,f}\,(u^n-f+\alpha\Delta f)^2. \]
This introduces a soft-thresholding with respect to the first approach. To sum up, we can choose the interpolation data such that the pixel density is increasing with $|u^n-f+\alpha\Delta f|$.
\section{The Iterative Methods}
\label{sec:algorithms}

With the model proposed in previous sections it is clear that we do need the complete sequence of optimal set $(K_n)_n$ in order to reconstruct the solutions $(u^n)_n$. While it does not represent a problem for denoising only purpose, it is not conceivable to store the complete sequence of inpainting mask for compression purpose. To overcome this issue, we propose in this section three algorithms, with on one side the hard-threshold criteria, namely we select the pixels where $|u^n-f+\alpha\Delta f|$ is maximum and on the other side, the soft-threshold criteria of the fat pixels approach, where the selected pixels are chosen according to the distribution of $|u^n-f+\alpha\Delta f|$.

\subsection{L2-INSTA}

\paragraph{Encoding.} During the encoding step, the original image is known on the whole domain $D$. We can thus set $u_0=f$ in $D$. We provide the input image $f$ and the desired pixel density $0<c<1$, and the algorithm produces the ``optimal'' mask $K$. We notice that $K$ is only optimal in the sense that it is a limit of optimal sets $K_n$. At each iteration, we compute and use a new mask without preserving the previous set. The stopping criterium is a given integer $N$ corresponding to the level of error on the reconstruction we want to achieve. We note that with this algorithm, $K_0$ is the same set as in the stationary case \cite{Belhachmi2022}. The complexity of the encoding is $O(N)$. We give the algorithm of \textit{L2-INSTA} in Algorithm \ref{algo:mask-creation:l2insta}.

\IncMargin{1.5em}
\begin{algorithm}[H]
    \SetAlgoLined
    \KwData{Original (noisy) image $f$, parameter (time-step) $\alpha>0$, desired pixel density $c\in (0,1)$, number of iterations $N\in \mathbb{N}^*$.}
    \KwResult{Inpainting mask $K\subset D$, last encoding reconstruction $u^N$.}
    
    $u^0\leftarrow f$\;
    
    \For{$n$ in $\{0,\hdots, N-1\}$}{
        Save in $K$ the $c\,|D|$ pixel by using the hard/soft-thresholding criterium i.e. $|u^n-f+\alpha\Delta f|$\;
        
        Compute $u^{n+1}$, solution of Problem \ref{pb:l2insta}\;
    }

    \caption{\textit{L2-INSTA}.}
    \label{algo:mask-creation:l2insta}
\end{algorithm}
\DecMargin{1.5em}

\paragraph{Decoding.} During the decoding step, the data are only available on $K$. Therefore, we must not use the image $f$ in $D\setminus K$ for the inpainting problem. We propose to put $u_0$ to zero in $D\setminus K$. However, doing this choice leads to a reconstruction $u$ which is near zero in the unknown set when computing the homogeneous heat equation for small time $t$. Since we want to compute the reconstruction at the same time as in the encoding step, we have to choose $t=N\alpha$. Thus, to have large $t$, we can either have a large number of iteration $N$, either a large time-step $\alpha$. But, for large $\alpha$, the criterium will be close to $|\Delta f|$ and the resulting mask will not depends on $u^n$ anymore. Thus, a good choice would be to have a small time-step $\alpha$ and a large iteration number $N$ in the encoding step. Note that for the decoding step, we can directly compute the reconstruction at time $t=N\alpha$ by performing only one iteration.

\subsection{L2-DEC}

\paragraph{Encoding.} In their article, Mainberger \textit{et al} proposed a probabilistic algorithm to compute an inpainting mask for a given inpainting method \cite{Mainberger2012}. This \textit{sparsification} algorithm is described as follow : we start with a mask containing all the pixels of the image, then at each iteration, we randomly delete a fraction of the pixels of the mask, we inpaint, we calculate the local error at each deleted pixel and we put back in the mask a subset of the deleted pixels presenting the most important error. We propose a modified algorithm based on the \textit{sparsification} one that we call \textit{L2-DEC}. The difference here is that, instead of selectionning a small random candidate set of pixels for each iteration, we directly choose an optimal set with respect to the criterion stated in this article for the current iteration. Using this strategy eliminate the probabilistic nature of the algorithm. The complexity of the encoding is $O\big(\frac{c}{q}\big)$, with $q$ the fraction of finally deleted pixel at the end of the iteration. Despite we give the algorithm of \textit{L2-DEC} in Algorithm \ref{algo:mask-creation:l2dec} for completness, we believe that \textit{L2-DEC} is not pertinent. Indeed, for every point $x_0$ in $K_n$, $u^n(x_0)=f(x_0)$ and the criterion becomes
\[ |u^n(x_0)-f(x_0)+\alpha\Delta f(x_0)| = |f(x_0)-f(x_0)+\alpha\Delta f(x_0)| = \alpha|\Delta f(x_0)|. \]
The algorithm will remove pixel with the lowest values $|\Delta f|$ and produce the same mask as the optimal one for the homogeneous diffusion inpainting \cite{Belhachmi2009}. We propose this theoretical algorithm for completeness purpose, but we will not realize numerical computations.

\IncMargin{1.5em}
\begin{algorithm}[H]
    \SetAlgoLined
    \KwData{Original (noisy) image $f$, parameter (time-step) $\alpha>0$, fraction $q$ of pixel removed, desired pixel density $c\in (0,1)$.}
    \KwResult{Inpainting mask $K\subset D$, last encoding reconstruction $u^N$.}
    
    $u^0\leftarrow f$\;
    $K\leftarrow D$\;
    $n\leftarrow 0$\;
    \While{$|K|>c\,|D|$}{
        Keep in $K$ the $|K|-q\,|D|$ pixel in $K$ by using the hard/soft-thresholding criterium i.e. $|u^n-f+\alpha\Delta f|$\;
        
        Compute $u^{n+1}$, solution of Problem \ref{pb:l2insta}\;
        
        $n\leftarrow n+1$\;
    }
    
    \caption{\textit{L2-DEC}.}
    \label{algo:mask-creation:l2dec}
\end{algorithm}
\DecMargin{1.5em}

\paragraph{Decoding.} As described in the original paper of the \textit{sparsification} algorithm \cite{Mainberger2012}, we should use the same inpainting (with the same parameters) as the one used in the encoding iterations, namely the discretized heat equation. However, as explained in the previous section, in order to have a pertinent reconstruction we should take the time-step $\alpha$ large enough. Such a choice will lead to a constant inpainting mask with respect to the iterations, which is not the purpose of this paper. We propose instead to use the homogeneous diffusion inpainting \cite{Belhachmi2009, Mainberger2012} to reconstruct the image from data in $K$.

\subsection{L2-INC}

\paragraph{Encoding.} Peter \textit{et al} proposed the \textit{densification} algorithm \cite{Peter2016Nov, Adam2017}. We start with an empty mask and, for each iteration, we randomly choose $\alpha$ pixels that do not belong to the mask. We then add only the single pixel that improves the global reconstruction error with respect to the image. In order to make the comparison possible with the following proposed algorithm, we propose to not only add one pixel but to add the fraction $q$ of pixels that improves the most the global reconstruction error with respect to the image. With our proposed algorithm, we directly choose an optimal candidate set with respect to the criterion stated in this article for the current iteration instead of choosing a small random candidate set of pixel for each iteration.  The complexity of the encoding is $O\big(\frac{c}{q}\big)$. We give the algorithm of \textit{L2-INC} in Algorithm \ref{algo:mask-creation:l2inc}.

\IncMargin{1.5em}
\begin{algorithm}[H]
    \SetAlgoLined
    \KwData{Original (noisy) image $f$, parameter (time-step) $\alpha>0$, fraction $q$ of pixel added, desired pixel density $c\in (0,1)$.}
    \KwResult{Inpainting mask $K\subset D$, last encoding reconstruction $u^N$.}
    
    $u^0\leftarrow f$\;
    $K\leftarrow \emptyset$\;
    $n\leftarrow 0$\;
    
    \While{$|K|<c\,|D|$}{
        Add the $q\,|D|$ pixels in $D\setminus K$ by using the hard/soft-thresholding criterium i.e. $|u^n-f+\alpha\Delta f|$, and add them to $K$\;
        
        
        Compute $u^{n+1}$, solution of Problem \ref{pb:l2insta}\;
        
        $n\leftarrow n+1$\;
    }
    
    \caption{\textit{L2-INC}.}
    \label{algo:mask-creation:l2inc}
\end{algorithm}
\DecMargin{1.5em}

\paragraph{Decoding.} With the same reasoning as the decoding step of the \textit{L2-DEC} algorithm, we propose to use the homogeneous diffusion inpainting to reconstruct the image from data in $K$. 
\section{Numerical Comparison of the Iterative Methods}
\label{sec:comparison}

In this section, we will present numerical results and comparisons between the presented methods from Section \ref{sec:algorithms}. The soft-thresholding rule from the ``fat pixel'' point of view can be enforced with a standard digital halftoning. Digital halftoning is a method of rendering that convert a continuous image to a binary image while giving the illusion of color continuity \cite{Ulichney1987, Adler2003}. This color continuity is simulated for the human eye by a spacial distribution of black and white pixels. An ideal digital halftoning method would preserves the average value of gray while giving the illusion of color continuity. For the experiment using the soft-thresholding rule, we use the Floyd-Steinberg dithering \cite{Floyd1976} except for \textit{L2-INC-} for reasons that will be discussed in the related section. We use the notation ``name-of-the-algorithm-\textit{T}'' and ``name-of-the-algorithm-\textit{H}'' to distinguish the hard- and soft-thresholding criteria used.

\subsection{Image Compression}

As discussed in the previous section, for \textit{L2-STA-} algorithms, we take for $\alpha=0.05$ (a small value). For \textit{L2-INC-H}, we cannot use the Floyd-Steinberg algorithm since it is biased by the error propagation when the number of pixel demanded is to small (here, we choose to add $50$ pixels per iteration). This leads to select pixel localised in the bottom-right part of the image. As an alternative, we propose to use a halftoning method based on the Lloyd’s method instead of error propagation \cite{Secord2002}. Algorithms \textit{H1-T} and \textit{H1-H} correspond to the methods described in \cite{Belhachmi2009}, namely, the hard/soft-thresholding of the absolute value of the laplacian of the input image $f_\delta$ i.e. $|\Delta f_\delta|$. These two methods have the advantage to not require any reconstructions during the encoding phase. However, since our input image may contain noise, the laplacian will be highly perturbed and the information given by the edges will be lost. It will lead to a poor reconstructions for images that contain noise.

In most of the case, \textit{L2-INC-T} gives greater or similar reconstruction quality than the other methods. It can be explained by the fact that, for each iteration, we add to the mask a small amount of best pixels for the current iteration. Thus, we keep important pixels for previous iterations and somehow store the full sequence of optimal sets $(K_n)_n$ containing a few number of pixels. Contrasting with them, the \textit{L2-INSTA-} methods ``forget'' at each iteration the optimal pixel set from the previous iterations. This leads to a poorer visual quality in the reconstruction. It is also interesting to notice that our model and then our analytical criteria are more robust to noise than the one using the homogeneous diffusion as inpainting operator, namely the \textit{H1-} methods. It confirm the pertinence of our model to handle image with gaussian noise.

In Table \ref{tab:methods-comparison:0.05} and Table \ref{tab:methods-comparison:0.10}, we give the $L^2$-errors between the image $f$ and the reconstruction $u$ (we write $L^2$ in the tables), as a function of the gaussian noise standard deviation for each method for different size of mask. In Figure \ref{fig:methods-comparison:0.1:wn:0}, Figure \ref{fig:methods-comparison:0.1:wn:0.03} and Figure \ref{fig:methods-comparison:0.1:wn:0.05} (Appendix \ref{app:images-compression}), we present various masks obtained and the corresponding reconstructed images.

\begin{table}[H]
    \centering
    \begin{tabular}{|c||c||c||c|c||c|c||c|c||c|c|}
        \hline
        \multirow{2}{*}{\textbf{Noise}} & \multicolumn{1}{c||}{\textbf{H1-T}} & \multicolumn{1}{c||}{\textbf{H1-H}} & \multicolumn{2}{c||}{\textbf{L2-INSTA-T}} & \multicolumn{2}{c||}{\textbf{L2-INSTA-H}} & \multicolumn{2}{c||}{\textbf{L2-INC-T}} & \multicolumn{2}{c|}{\textbf{L2-INC-H}} \\
        \cline{2-11}
              & $L^2$ & $L^2$ & $N$ & $L^2$ & $N$ & $L^2$ & $\alpha$ & $L^2$ & $\alpha$ & $L^2$ \\
        \hhline{|===========|}
           0  & 34.46 & \textbf{11.43} & 3131 & 29.39 & 6517 & 11.58 & 0.26 & 18.68 & 2.25 & 20.18 \\
        \hline
         0.03 & 15.84 & 14.21 & 7905 & 15.76 & 8765 & 13.37 & 0.26 &  \textbf{8.46} & 0.42 & 21.24 \\
        \hline
         0.05 & 19.24 & 16.44 & 4341 & 18.82 & 6472 & 15.50 & 0.16 & \textbf{11.74} & 1.44 & 24.25 \\
        \hline
         0.1  & 30.91 & 23.78 & 4507 & 24.81 & 6192 & \textbf{21.52} & 0.06 & 23.27 & 1.43 & 27.22 \\
        \hline
         0.2  & 64.51 & 40.68 & 4717 & 40.24 & 6376 & 35.26 & 0.11 & 48.31 & 0.62 & \textbf{33.06} \\
        \hline
    \end{tabular} \vspace{0.2cm}

    \caption{$L^2$-error between the original image $f$ and the reconstruction $u$ (build from  $f_\delta$) with $5\%$ of total pixels saved.}
    \label{tab:methods-comparison:0.05}
\end{table}

\begin{table}[H]
    \centering
    \begin{tabular}{|c||c||c||c|c||c|c||c|c||c|c|}
        \hline
        \multirow{2}{*}{\textbf{Noise}} & \multicolumn{1}{c||}{\textbf{H1-T}} & \multicolumn{1}{c||}{\textbf{H1-H}} & \multicolumn{2}{c||}{\textbf{L2-INSTA-T}} & \multicolumn{2}{c||}{\textbf{L2-INSTA-H}} & \multicolumn{2}{c||}{\textbf{L2-INC-T}} & \multicolumn{2}{c|}{\textbf{L2-INC-H}} \\
        \cline{2-11}
              & $L^2$ & $L^2$ & $N$ & $L^2$ & $N$ & $L^2$ & $\alpha$ & $L^2$ & $\alpha$ & $L^2$ \\
        \hhline{|===========|}
           0  & 22.99 &  \textbf{6.02} & 3655 & 13.61 & 4815 &  6.78 & 0.31 &  6.36 & 1.03 & 14.64 \\
        \hline
         0.03 & 11.52 & 10.43 & 2315 & 13.69 & 2407 & 10.49 & 0.11 &  \textbf{6.96} & 0.42 & 17.50 \\
        \hline
         0.05 & 15.42 & 13.77 & 3131 & 14.50 & 3909 & 12.64 & 0.11 & \textbf{10.94} & 0.83 & 20.89 \\
        \hline
         0.1  & 27.51 & 23.09 & 2143 & 22.17 & 5060 & \textbf{20.34} & 0.11 & 20.72 & 0.83 & 23.56 \\
        \hline
         0.2  & 55.09 & 42.44 & 3007 & 35.60 & 4987 & 35.99 & 0.11 & 39.14 & 0.83 & \textbf{31.87} \\
        \hline
    \end{tabular} \vspace{0.2cm}

    \caption{$L^2$-error between the original image $f$ and the reconstruction $u$ (build from  $f_\delta$) with $10\%$ of total pixels saved.}
    \label{tab:methods-comparison:0.10}
\end{table}

\begin{figure}[H]
	\centering
	\subfloat[Edges with \textit{L2-INC-T} method.]{
		\includegraphics[height=3.6cm]{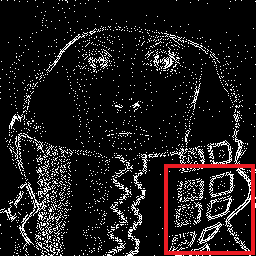}
	}
	\quad
	\subfloat[Zoom.]{
		\includegraphics[height=3.6cm]{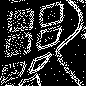}
	}
	\quad
	\subfloat[Edges with \textit{L2-STA-H} method.]{
		\includegraphics[height=3.6cm]{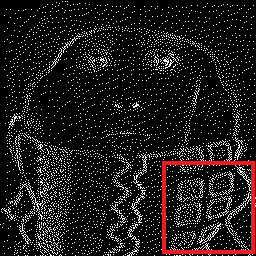}
	}
	\quad
	\subfloat[Zoom.]{
		\includegraphics[height=3.6cm]{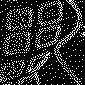}
	}
	\caption{Zoom on homogeneous parts and edges.}
	\label{fig:zoom}
\end{figure}

In Figure \ref{fig:zoom} we see in particular that, in order to reconstruct an homogeneous part of the image, \textit{L2-INC-T} select a few amount of pixels (here $4$ pixels) whereas an halftoning algorithm way more redundant pixels. It leads to more pixels available to spend near the edges in \textit{L2-INC-T} and then, a better reconstruction.

\begin{figure}[H]
	\centering
	\subfloat[Homogeneous part with \textit{L2-INC-T} method.]{
		\includegraphics[height=3.6cm]{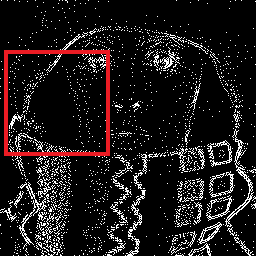}
	}
	\quad
	\subfloat[Zoom.]{
		\includegraphics[height=3.6cm]{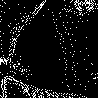}
	}
	\quad
	\subfloat[Homogeneous part with \textit{L2-STA-H} method.]{
		\includegraphics[height=3.6cm]{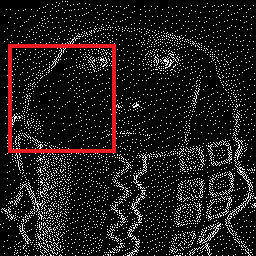}
	}
	\quad
	\subfloat[Zoom.]{
		\includegraphics[height=3.6cm]{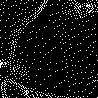}
	}
\end{figure}

\subsection{Image Denoising}

We start by making the hereafter observation : for well chosen parameters, the error between the original image $f$ and the solution $u^n$ decreases during the encoding step and becomes lower than the error between the original image $f$ and the noised one $f_\delta$. In fact, this error on $u^n$ is lower than the one on $f_\delta$ from the first iteration. We propose then to use the last encoding reconstruction $u^N$ of algorithms \textit{L2-INSTA} and \textit{L2-INC}, Algorithm \ref{algo:mask-creation:l2insta} and Algorithm \ref{algo:mask-creation:l2inc} respectively, as denoised version of $f_\delta$. For every noise level $\sigma$, we take take $1\%$ of the total pixel in the mask for \textit{L2-INSTA-} and $2\%$ for \textit{L2-INC-} methods except when $\sigma=0.2$, we take $4\%$ for the \textit{L2-INC-} methods. In any case, we set $\alpha=0.01$.

It appears that \textit{L2-INC-} are almost as efficient as the linear diffusion filter, which gives the lowest error, but are more edges preserving than the linear diffusion filter. As expected, every of the proposed methods perform noise removing and we will exploit this feature in the next section.

We give in Table \ref{tab:denoising-comparison} the $L^2$-error between $f$ and $u^N$ with respect to noise level in $f_\delta$ for the methods proposed in this paper and for the linear diffusion filter of parameter $\eta$, which is known to remove gaussian noise, and the resulting images in Figure \ref{fig:denoising} (Appendix \ref{app:images-denoising}).

\begin{table}[H]
    \centering
    \begin{tabular}{|c||c||c|c||c|c||c|c||c||c|}
        \hline
        \multirow{2}{*}{\textbf{Noise}} & \multirow{2}{*}{$\|f-f_\delta\|_{L^2(D)}$} & \multicolumn{2}{c||}{\textbf{Lin. Filter}} & \multicolumn{2}{c||}{\textbf{L2-INSTA-T}} & \multicolumn{2}{c||}{\textbf{L2-INSTA-H}} & \multicolumn{1}{c||}{\textbf{L2-INC-T}} & \multicolumn{1}{c|}{\textbf{L2-INC-H}} \\
        \cline{3-10}
             &  & $\eta$ & $L^2$ & $N$ & $L^2$ & $N$ & $L^2$ & $L^2$ & $L^2$ \\
        \hhline{|==========|}
         0.03 &  7.65 & 0.9 &  4.61 & 44 &  5.28 &  38 &  4.92 &  4.82 &  \textbf{4.44} \\
        \hline
         0.05 & 12.81 & 1.2 & \textbf{6.16} & 58 &  8.54 &  55 &  7.58 &  7.60 &  6.60 \\
        \hline
         0.1  & 25.45 & 2.0 & \textbf{8.98} & 59 & 16.60 & 102 & 14.06 & 14.67 & 12.50 \\
        \hline
         0.2  & 48.43 & 2.5 & \textbf{12.72} & 88 & 30.65 & 182 & 25.76 & 19.92 & 16.44 \\
        \hline
    \end{tabular} \vspace{0.2cm}
    \caption{Using $u^N$ as a denoised version of the input image.}
    \label{tab:denoising-comparison}
\end{table}

\section{Numerical Comparison with the ``Probabilistic'' Methods}
\label{sec:numerical_results}

In this section, we present numerical results and comparisons with some stat-of-the-art compression by inpainting methods for image with noise. In the following we denote by \textit{SPAR} the \textit{sparsification} algorithm without the \textit{nonlocal pixel exchange} post-optimization (step within every \textit{sparsification}'s iteration to avoid source locality due to local error computation) from \cite{Mainberger2012} and by \textit{DENS} a modified version of the original \textit{densification} algorithm from \cite{Adam2017}, that we will describe in the sequel.

The knowledge of the previous section motives us to replace $f$ by $u^n$ as Dirichlet boundary condition in the inpainting mask $K$, and thus, to replace $f$ by $u^n$ in the shape optimization problem \eqref{pb:l2insta_opt_no_constraint} : for a given $n\in\N$,
\begin{equation}
	\min_{K_n\subseteq D,\ \capop(K_n)\leq c}\Big\{  \frac{1}{2}\int_D |u_{K_n}-u^n|^2\ dx + \frac{\alpha}{2}\int_D |\nabla (u_{K_n}-u^n)|^2\ dx\ \Big\},
	\label{pb:l2insta_opt_no_constraint_modified}
\end{equation}
with $u_{K_n}$ solution of \\

\begin{problem} For $n\in\N$, given $u^n$, find $u^{n+1}$ in $H^1(D)$ such that
	\begin{equation}
		\left\{\begin{array}{rl}
			u^{n+1} - \alpha \Delta u^{n+1} = u^n, & \text{in}\ D\setminus K_n, \\
			u^{n+1} = u^n, & \text{in}\ K_n, \\
			\frac{\partial u^{n+1}}{\partial \mathbf{n}} = 0, & \text{on}\ \partial D, \\
		\end{array}\right .
	\end{equation}
\end{problem}

Since the analysis remains the same as in Section \ref{sec:l2insta:continuous_model}, these changes yield to a new criterion, namely $|\Delta u^n|$ for both hard- and soft-thresholding. This \textit{tonal optimization} avoid brutal smoothing of the image during the first iterations. \textit{Tonal optimization} consists in changing the value of pixel in a given inpainting mask in order to improve the overall image's reconstruction quality \cite{Mainberger2012}. Like for the previous algorithms, we take $u_0=f$ in $D$ for the encoding step since the entire noisy image is available, but we set $u_0$ to zero on $D\setminus K$ during the decoding step.

As suggested by the authors in the original paper of the \textit{SPAR} algorithm, the candidate set have $2\%$ of the pixel in $D$, and we choose the finally removed numbers of pixel to be $50$ pixels from the candidate set. To make the comparison possible, we propose to remove/add a fixed number of pixel (we choose $50$ pixels) at each iteration for the appropriate algorithms. Then, we have to modify the \textit{DENS} algorithm, since it originally add one single pixel per iteration. This modification induces lower visual quality for the reconstruction with a high gain for the computation time. Also, we choose the number of randomly selected candidates per iteration to be $100$ which still require a lot of computation time. Therefore, we have $100$ reconstructions to compute at each iteration. We name the method from this new model combined with Algorithm \ref{algo:mask-creation:l2inc} and hard-thresholding \textit{L2-INC-T-E}. 

\textit{L2-INC-T-E} seems to efficiently denoise the input image while compressing and outperform the \textit{SPAR} method and the \textit{DENS} method for high level of noise. However, our new model is designed to handle input image with gaussian noise and then, perform worst than the previous model when the image does not contain noise.

We give in Table \ref{tab:denoising-comparison} the $L^2$-error between $f$ and $u^N$ with respect to noise level in $f_\delta$ for the \textit{L2-INC-T-E} method proposed in this section and the ``probabilistic'' methods from the state-of-the-art and we give in Figure \ref{fig:probabilistic} and Figure \ref{fig:l2-inc-t-e} (Appendix \ref{app:denoising-compression}) some illustration. 



\begin{table}[H]
    \centering
    \begin{tabular}{|c||c||c||c|c|}
        \hline
        \multirow{2}{*}{\textbf{Noise}} & \multicolumn{1}{c||}{\textbf{SPAR}} & \multicolumn{1}{c||}{\textbf{DENS}} & \multicolumn{2}{c|}{\textbf{L2-INC-T}} \\
        \cline{2-5}
              & $L^2$ & $L^2$ & $\alpha$ & $L^2$ \\
        \hhline{|=====|}
           0  &  \textbf{3.68} &  8.70 & 0.01 & 15.54 \\
        \hline
         0.03 &  \textbf{6.94} &  9.09 & 0.01 & 7.53 \\
        \hline
         0.05 & 11.11 &  9.82 & 0.01 & \textbf{7.98} \\
        \hline
         0.1  & 21.37 & 12.58 & 0.01 & \textbf{11.79} \\
        \hline
    \end{tabular} \vspace{0.2cm}

    \caption{$L^2$-error between the original image $f$ and the reconstruction $u$ (build from  $f_\delta$ and with the new model) with $10\%$ of total pixels saved.}
    \label{tab:enhanced-methods-comparison:0.10}
\end{table}

\section*{Summary and Conclusions}
\label{sec:l2insta:conclusion}

In this article, we have considered a new fully parabolic model based on the heat equation for PDE images compression and formulate an associated shape optimization problem for the best choice of data interpolation choice. By using the results found in \cite{Belhachmi2022}, we proved that it has a relaxed solution in the framework of the $\gamma$-convergence and we proposed two points of view for choosing an optimal inpainting mask : a first one that use an asymptotic development similar to the topological gradient, and a second one by considering the sets formed by a finite number of balls that we called ``fat pixels''. Unlike the stationary case, the iterative approach include the results of the previous steps in the analytic criterion giving the best choice which is now depending on the iterations. Next, we proposed and implemented three algorithms to construct the sequences of inpainting mask $(K_n)_n$ : \textit{L2-INSTA}, \textit{L2-DEC} and \textit{L2-INC}, and we performed several simulations with  different compression rates and different amounts of noise for grayscale images. It appears that in most cases, masks issued form \textit{L2-INC-T} outperforms the other ones. Finally, we performed a tonal optimization within the same method by changing the Dirichlet boundary conditions on the mask, which appears to give better results for denoising while compressing an image with noise than the \textit{sparsification} and the \textit{densification} algorithms \cite{Adam2017, Mainberger2012} for high level of noise.

\bibliographystyle{siam}  
\bibliography{template}

\appendix

\newpage\section{Inpainting Masks and Reconstructions for Compression from Section \ref{sec:comparison}}
\label{app:images-compression}

\begin{figure}[H]
	\centering
	\subfloat[Mask with \textit{H1-T} method.]{
		\includegraphics[height=3.6cm]{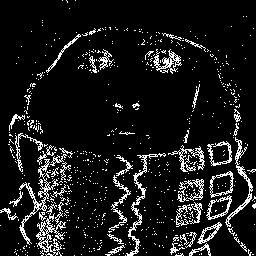}
	}
	\quad
	\subfloat[Reconstruction with \textit{H1-T} method.]{
		\includegraphics[height=3.6cm]{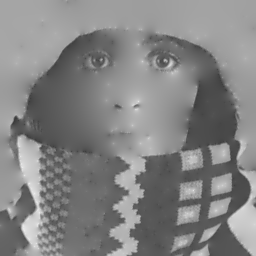}
	}
	\quad
	\subfloat[Mask with \textit{H1-H} method.]{
		\includegraphics[height=3.6cm]{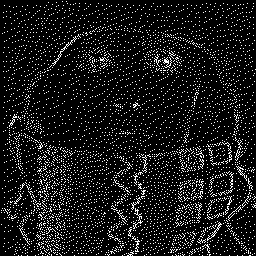}
	}
	\quad
	\subfloat[Reconstruction with \textit{H1-H} method.]{
		\includegraphics[height=3.6cm]{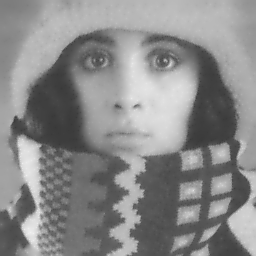}
	}
	\quad
	\subfloat[Mask with \textit{L2-INSTA-T} method.]{
		\includegraphics[height=3.6cm]{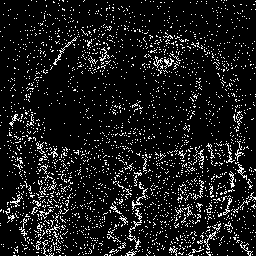}
	}
	\quad
	\subfloat[Reconstruction with \textit{L2-INSTA-T} method.]{
		\includegraphics[height=3.6cm]{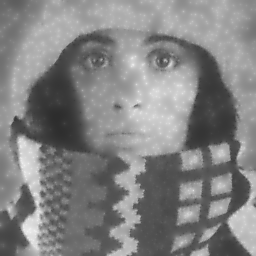}
	}
	\quad
	\subfloat[Mask with \textit{L2-INSTA-H} method.]{
		\includegraphics[height=3.6cm]{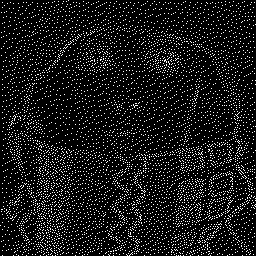}
	}
	\quad
	\subfloat[Reconstruction with \textit{L2-INSTA-H} method.]{
		\includegraphics[height=3.6cm]{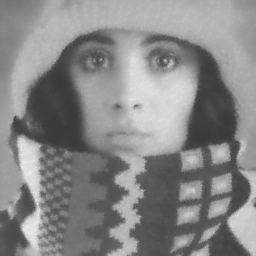}
	}
	\quad
	\subfloat[Mask with \textit{L2-INC-T} method.]{
		\includegraphics[height=3.6cm]{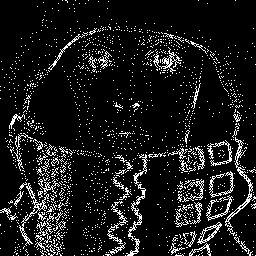}
	}
	\quad
	\subfloat[Reconstruction with \textit{L2-INC-T} method.]{
		\includegraphics[height=3.6cm]{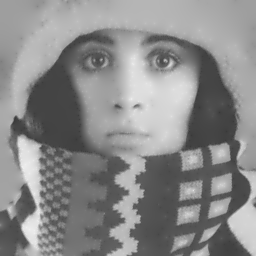}
	}
	\quad
	\subfloat[Mask with \textit{L2-INC-H} method.]{
		\includegraphics[height=3.6cm]{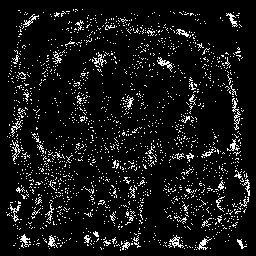}
	}
	\quad
	\subfloat[Reconstruction with \textit{L2-INC-H} method.]{
		\includegraphics[height=3.6cm]{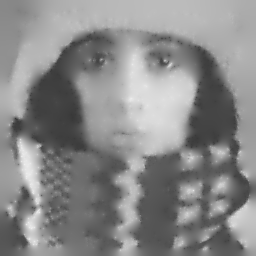}
	}
	
	\caption{Masks and reconstructions for Table \ref{tab:methods-comparison:0.10} when the input image is noiseless ($\sigma=0$).}
	\label{fig:methods-comparison:0.1:wn:0}
\end{figure}

\vfill\newpage

\begin{figure}[H]
	\centering
	\subfloat[Mask with \textit{H1-T} method.]{
		\includegraphics[height=3.6cm]{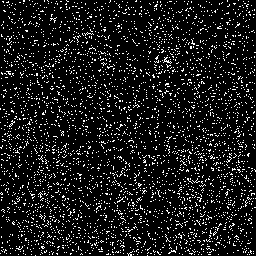}
	}
	\quad
	\subfloat[Reconstruction with \textit{H1-T} method.]{
		\includegraphics[height=3.6cm]{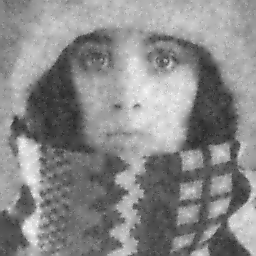}
	}
	\quad
	\subfloat[Mask with \textit{H1-H} method.]{
		\includegraphics[height=3.6cm]{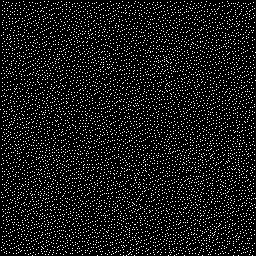}
	}
	\quad
	\subfloat[Reconstruction with \textit{H1-H} method.]{
		\includegraphics[height=3.6cm]{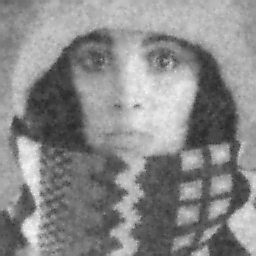}
	}
	\quad
	\subfloat[Mask with \textit{L2-INSTA-T} method.]{
		\includegraphics[height=3.6cm]{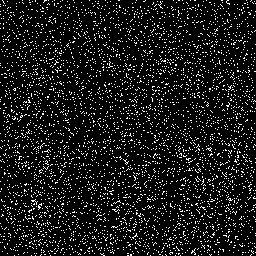}
	}
	\quad
	\subfloat[Reconstruction with \textit{L2-INSTA-T} method.]{
		\includegraphics[height=3.6cm]{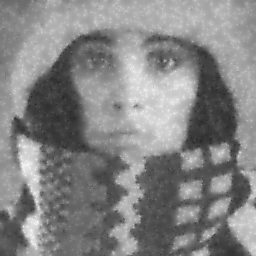}
	}
	\quad
	\subfloat[Mask with \textit{L2-INSTA-H} method.]{
		\includegraphics[height=3.6cm]{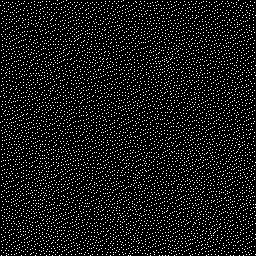}
	}
	\quad
	\subfloat[Reconstruction with \textit{L2-INSTA-H} method.]{
		\includegraphics[height=3.6cm]{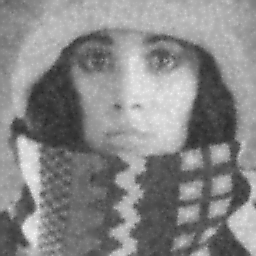}
	}
	\quad
	\subfloat[Mask with \textit{L2-INC-T} method.]{
		\includegraphics[height=3.6cm]{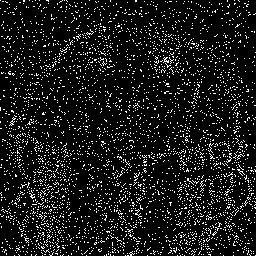}
	}
	\quad
	\subfloat[Reconstruction with \textit{L2-INC-T} method.]{
		\includegraphics[height=3.6cm]{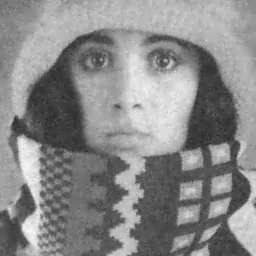}
	}
	\quad
	\subfloat[Mask with \textit{L2-INC-H} method.]{
		\includegraphics[height=3.6cm]{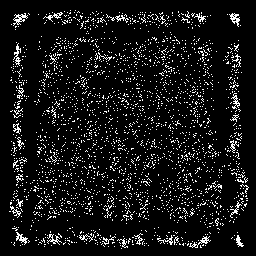}
	}
	\quad
	\subfloat[Reconstruction with \textit{L2-INC-H} method.]{
		\includegraphics[height=3.6cm]{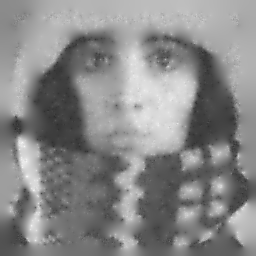}
	}
	
	\caption{Masks and reconstructions for Table \ref{tab:methods-comparison:0.10} when the input image is affected by gaussian noise ($\sigma=0.03$).}
	\label{fig:methods-comparison:0.1:wn:0.03}
\end{figure}

\begin{figure}[H]
	\centering
	\subfloat[Mask with \textit{H1-T} method.]{
		\includegraphics[height=3.6cm]{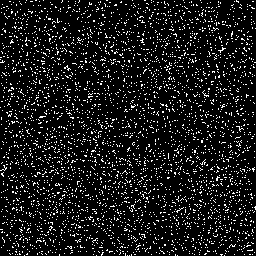}
	}
	\quad
	\subfloat[Reconstruction with \textit{H1-T} method.]{
		\includegraphics[height=3.6cm]{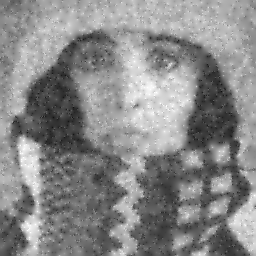}
	}
	\quad
	\subfloat[Mask with \textit{H1-H} method.]{
		\includegraphics[height=3.6cm]{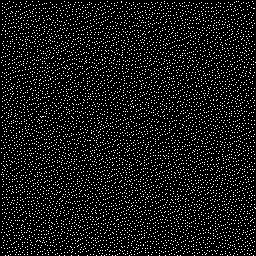}
	}
	\quad
	\subfloat[Reconstruction with \textit{H1-H} method.]{
		\includegraphics[height=3.6cm]{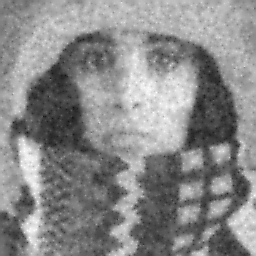}
	}
	\quad
	\subfloat[Mask with \textit{L2-INSTA-T} method.]{
		\includegraphics[height=3.6cm]{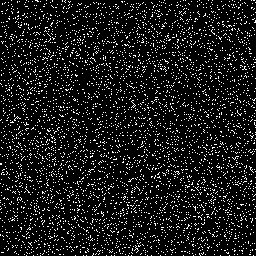}
	}
	\quad
	\subfloat[Reconstruction with \textit{L2-INSTA-T} method.]{
		\includegraphics[height=3.6cm]{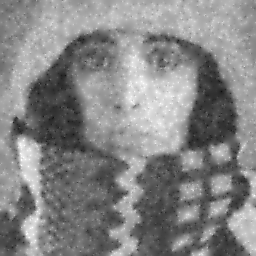}
	}
	\quad
	\subfloat[Mask with \textit{L2-INSTA-H} method.]{
		\includegraphics[height=3.6cm]{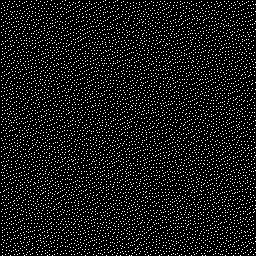}
	}
	\quad
	\subfloat[Reconstruction with \textit{L2-INSTA-H} method.]{
		\includegraphics[height=3.6cm]{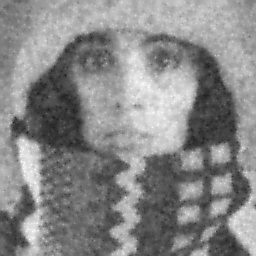}
	}
	\quad
	\subfloat[Mask with \textit{L2-INC-T} method.]{
		\includegraphics[height=3.6cm]{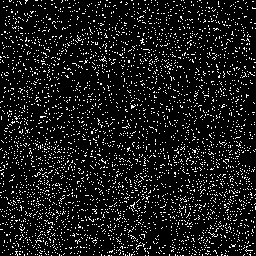}
	}
	\quad
	\subfloat[Reconstruction with \textit{L2-INC-T} method.]{
		\includegraphics[height=3.6cm]{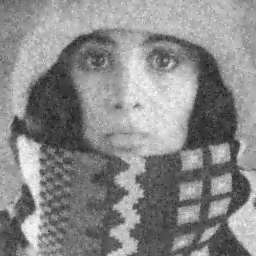}
	}
	\quad
	\subfloat[Mask with \textit{L2-INC-H} method.]{
		\includegraphics[height=3.6cm]{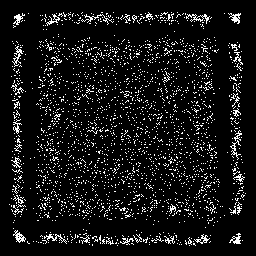}
	}
	\quad
	\subfloat[Reconstruction with \textit{L2-INC-H} method.]{
		\includegraphics[height=3.6cm]{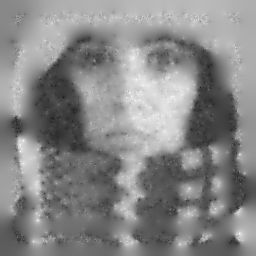}
	}
	
	\caption{Masks and reconstructions for Table \ref{tab:methods-comparison:0.10} when the input image is affected by gaussian noise ($\sigma=0.05$).}
	\label{fig:methods-comparison:0.1:wn:0.05}
\end{figure}
\newpage\section{Reconstructions for Image Denoising from Section \ref{sec:comparison}}
\label{app:images-denoising}

\begin{figure}[H]
	\centering
	\subfloat[\textit{L2-INSTA-T} ($\sigma=0.03$).]{
		\includegraphics[height=3.6cm]{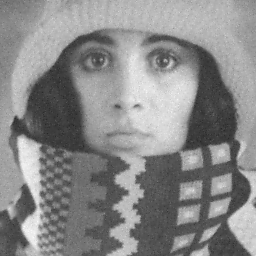}
	}
	\quad
	\subfloat[\textit{L2-INSTA-T} ($\sigma=0.05$).]{
		\includegraphics[height=3.6cm]{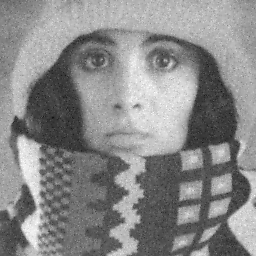}
	}
	\quad
	\subfloat[\textit{L2-INSTA-T} ($\sigma=0.1$).]{
		\includegraphics[height=3.6cm]{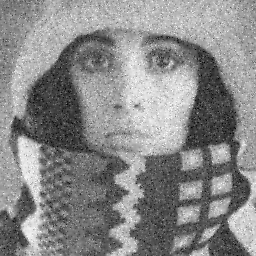}
	}
	\quad
	\subfloat[\textit{L2-INSTA-T} ($\sigma=0.2$).]{
		\includegraphics[height=3.6cm]{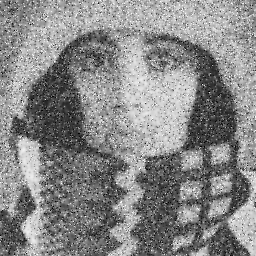}
	}
	\quad
	\subfloat[\textit{L2-INSTA-H} ($\sigma=0.03$).]{
		\includegraphics[height=3.6cm]{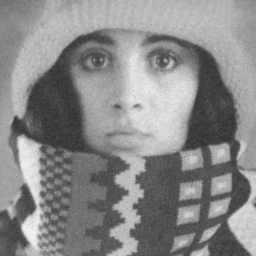}
	}
	\quad
	\subfloat[\textit{L2-INSTA-H} ($\sigma=0.05$).]{
		\includegraphics[height=3.6cm]{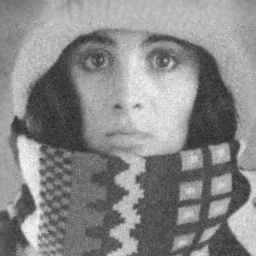}
	}
	\quad
	\subfloat[\textit{L2-INSTA-H} ($\sigma=0.1$).]{
		\includegraphics[height=3.6cm]{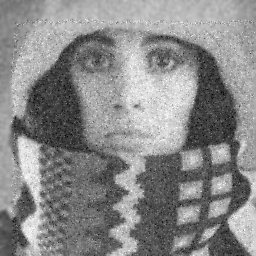}
	}
	\quad
	\subfloat[\textit{L2-INSTA-H} ($\sigma=0.2$).]{
		\includegraphics[height=3.6cm]{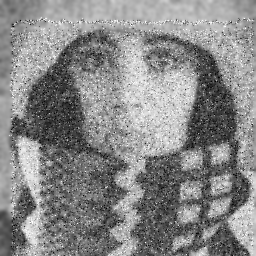}
	}
	\quad
	\subfloat[\textit{L2-INC-T} ($\sigma=0.03$).]{
		\includegraphics[height=3.6cm]{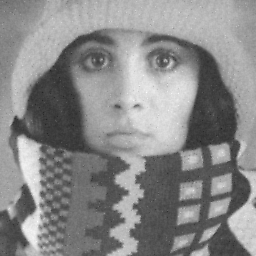}
	}
	\quad
	\subfloat[\textit{L2-INC-T} ($\sigma=0.05$).]{
		\includegraphics[height=3.6cm]{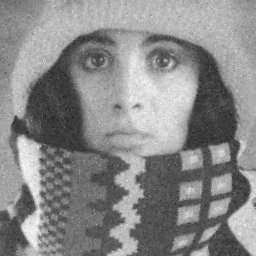}
	}
	\quad
	\subfloat[\textit{L2-INC-T} ($\sigma=0.1$).]{
		\includegraphics[height=3.6cm]{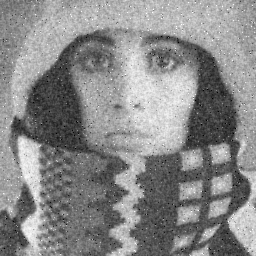}
	}
	\quad
	\subfloat[\textit{L2-INC-T} ($\sigma=0.2$).]{
		\includegraphics[height=3.6cm]{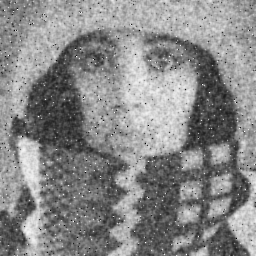}
	}
	\quad
	\subfloat[\textit{L2-INC-H} ($\sigma=0.03$).]{
		\includegraphics[height=3.6cm]{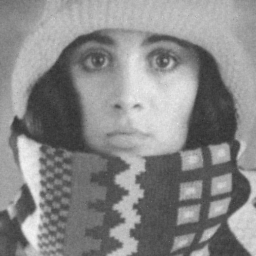}
	}
	\quad
	\subfloat[\textit{L2-INC-H} ($\sigma=0.05$).]{
		\includegraphics[height=3.6cm]{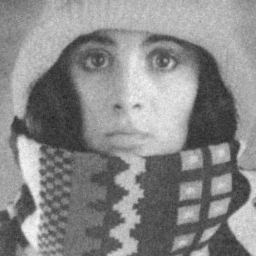}
	}
	\quad
	\subfloat[\textit{L2-INC-H} ($\sigma=0.1$).]{
		\includegraphics[height=3.6cm]{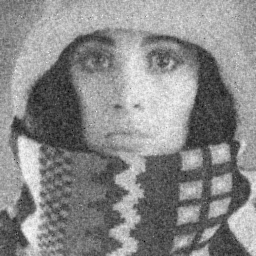}
	}
	\quad
	\subfloat[\textit{L2-INC-H} ($\sigma=0.2$).]{
		\includegraphics[height=3.6cm]{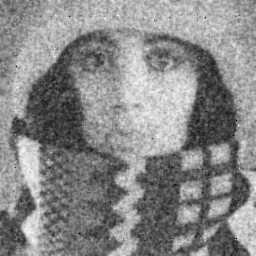}
	}
	\caption{$u^N$ as a denoised version of the input image for multiple levels of noise.}
	\label{fig:denoising}
\end{figure}
\newpage\section{Inpainting Masks and Reconstructions for the new model from Section \ref{sec:numerical_results}}
\label{app:denoising-compression}

\begin{figure}[H]
    \subfloat[Input ($\sigma=0)$.]{
		\includegraphics[width=3.6cm]{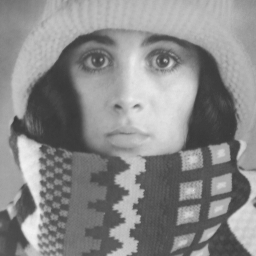}
	}
	\quad
	\subfloat[Input ($\sigma=0.03)$.]{
		\includegraphics[width=3.6cm]{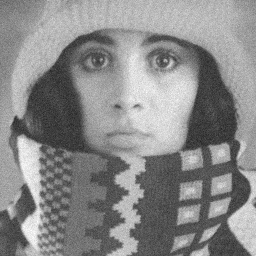}
	}
	\quad
	\subfloat[Input ($\sigma=0.05)$.]{
		\includegraphics[width=3.6cm]{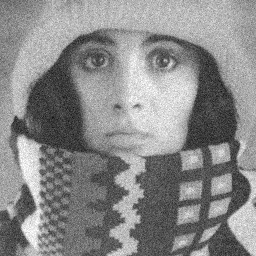}
	}
	\quad
	\subfloat[Input ($\sigma=0.1)$.]{
		\includegraphics[width=3.6cm]{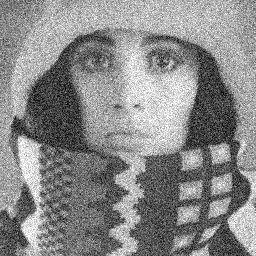}
	}
	\quad
	\subfloat[\textit{DENS}.]{
		\includegraphics[width=3.6cm]{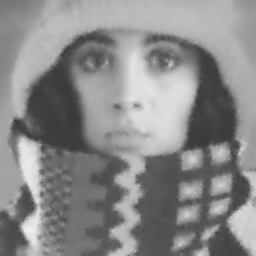}
	}
	\quad
	\subfloat[\textit{DENS}.]{
		\includegraphics[width=3.6cm]{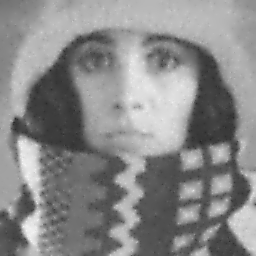}
	}
	\quad
	\subfloat[\textit{DENS}.]{
		\includegraphics[width=3.6cm]{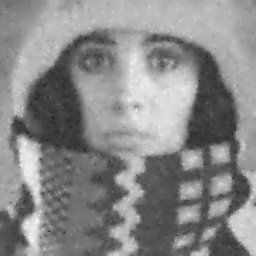}
	}
	\quad
	\subfloat[\textit{DENS}.]{
		\includegraphics[width=3.6cm]{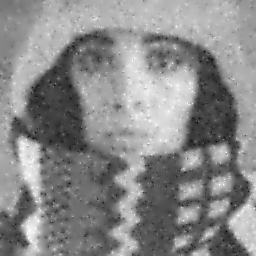}
	}
	\quad
	\subfloat[\textit{SPAR}.]{
		\includegraphics[width=3.6cm]{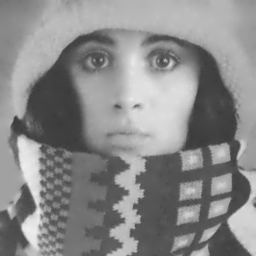}
	}
	\quad
	\subfloat[\textit{SPAR}.]{
		\includegraphics[width=3.6cm]{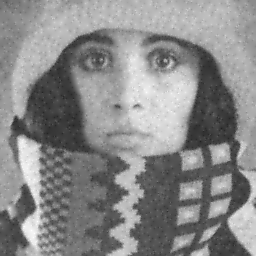}
	}
	\quad
	\subfloat[\textit{SPAR}.]{
		\includegraphics[width=3.6cm]{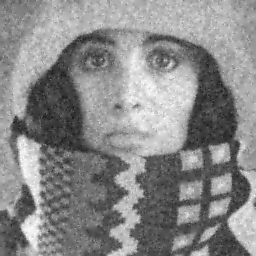}
	}
	\quad
	\subfloat[\textit{SPAR}.]{
		\includegraphics[width=3.6cm]{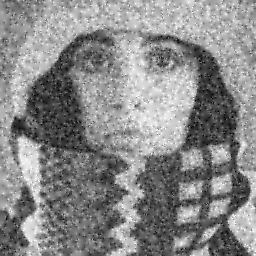}
	}
	\caption{Reconstruction with the \textit{sparsification} and \textit{densification} methods with $10\%$ of total pixels saved.}
	\label{fig:probabilistic}
\end{figure}

\vfill\newpage

\begin{figure}[H]
    \subfloat[Input ($\sigma=0)$.]{
		\includegraphics[width=3.6cm]{resources/images/Trui.png}
	}
	\quad
	\subfloat[Input ($\sigma=0.03)$.]{
		\includegraphics[width=3.6cm]{resources/images/Trui-wn-0.03.png}
	}
	\quad
	\subfloat[Input ($\sigma=0.05)$.]{
		\includegraphics[width=3.6cm]{resources/images/Trui-wn-0.05.png}
	}
	\quad
	\subfloat[Input ($\sigma=0.1)$.]{
		\includegraphics[width=3.6cm]{resources/images/Trui-wn-0.1.png}
	}
	\quad
	\subfloat[Mask.]{
		\includegraphics[width=3.6cm]{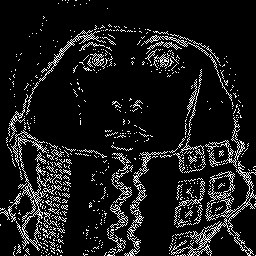}
	}
	\quad
	\subfloat[Mask.]{
		\includegraphics[width=3.6cm]{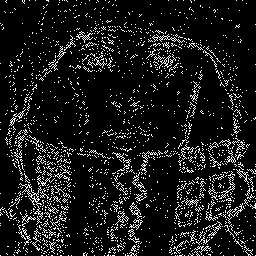}
	}
	\quad
	\subfloat[Mask.]{
		\includegraphics[width=3.6cm]{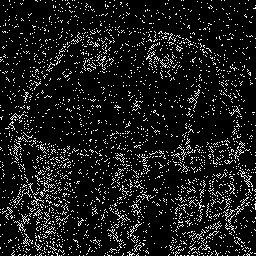}
	}
	\quad
	\subfloat[Mask.]{
		\includegraphics[width=3.6cm]{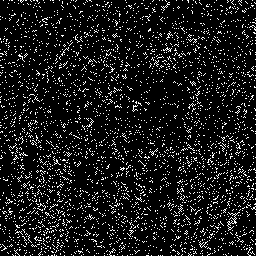}
	}
	\quad
	\subfloat[Data $u^N$.]{
		\includegraphics[width=3.6cm]{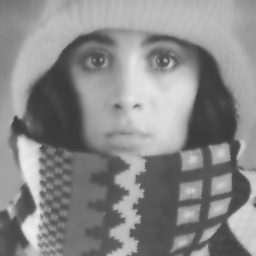}
	}
	\quad
	\subfloat[Data $u^N$.]{
		\includegraphics[width=3.6cm]{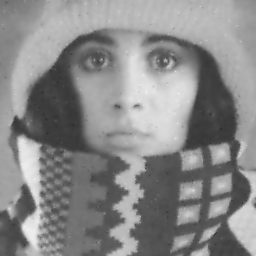}
	}
	\quad
	\subfloat[Data $u^N$.]{
		\includegraphics[width=3.6cm]{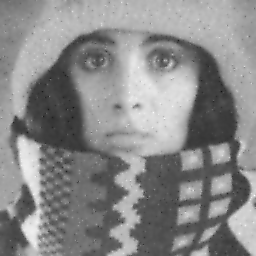}
	}
	\quad
	\subfloat[Data $u^N$.]{
		\includegraphics[width=3.6cm]{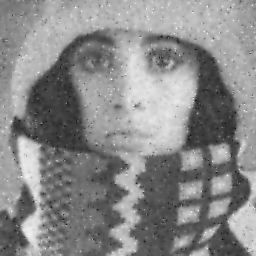}
	}
	\quad
	\subfloat[Reconstruction.]{
		\includegraphics[width=3.6cm]{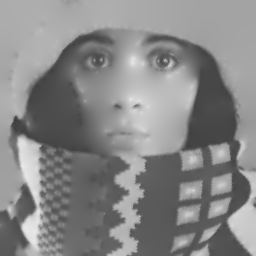}
	}
	\quad
	\subfloat[Reconstruction.]{
		\includegraphics[width=3.6cm]{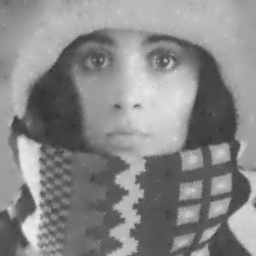}
	}
	\quad
	\subfloat[Reconstruction.]{
		\includegraphics[width=3.6cm]{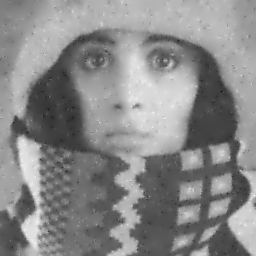}
	}
	\quad
	\subfloat[Reconstruction.]{
		\includegraphics[width=3.6cm]{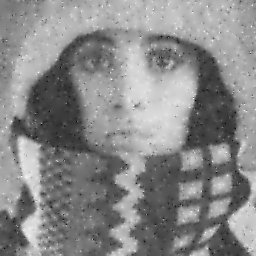}
	}
	\caption{Reconstruction with \textit{L2-INC-T-E} method with $10\%$ of total pixels saved.}
	\label{fig:l2-inc-t-e}
\end{figure}

\end{document}